\documentclass[12pt]{article}

\usepackage{latexsym,amssymb}
\usepackage{graphicx}
\usepackage{graphics}
\usepackage[all]{xy}
\CompileMatrices

\newcommand{\mz}{\mathbb Z}

\newcommand{\low}{\mbox{\it low}}

\def\dem{\noindent{\it Proof: \ }}
\def\QED{\hfill$\Box$}

\newtheorem{Teo}{Theorem}[section]
\newtheorem{Lema}[Teo]{Lemma}
\newtheorem{Coro}[Teo]{Corollary}
\newtheorem{Propo}[Teo]{Proposition}
\newtheorem{Remark}[Teo]{Remark}

\begin{document}

\title{\Large\bf
The $ku$-homology of certain classifying spaces $II$.
\vspace{-2mm}}

\author{\normalsize\sc
\thanks{Partially supported by Cimat and Fundaci\'on Koval\'evskaia.}
Leticia Z\'arate }

\date{}

\maketitle

\begin{abstract}
We calculate the annihilator of the $ku$-toral class for the
$p$-groups ${\mz}_{p^2} \times {\mz}_{p^k}$ with $k \geq 3$. This allows
us to give a description of the $ku_*$-homology of the groups we are
dealing with.
\end{abstract}

{\small


\noindent {{\it Proposed running head}: { On the $\mz_{p^k}
\times \mz_{p^e}$ Conner-Floyd conjecture.}} }

\section{Introduction.}\label{s1}

We are interested in the structure of complex bordism of the $p$-groups $\mz_{p^e} \times \mz_{p^k}$.
By Quillen's splitting theorem we know that it is enough to calculate the $BP$-homology.
Using the results in \cite{Nakos} we proved in \cite{kuosa} that for the groups $\mz_{p^2} \times \mz_{p^2}$
and $\mz_{p^e} \times \mz_{p}$ the $ku$-homology contains all the complex
bordism information. Indeed we constructed a set of generators
of the annihilator of the $ku$-toral class that are elements of $BP \langle  1 \rangle$ and that
also are a set of generators of the annihilator of the $BP$-toral class.

We conjectured in \cite{kuosa} that, for the groups of the form $\mz_{p^e} \times \mz_{p^k}$,
the annihilator of the $ku$-toral class gives all the complex bordism information. In this work we
obtain a set of generators for the annhilator of the $ku$-toral class for the groups $\mz_{p^2} \times \mz_{p^k}$
(with $k \geq 3$).

We obtain that for the group $\mz_{p^2} \times \mz_{p^{k+1}}$ the annihilator of the $ku$-toral class is given by:
\begin{equation}\label{an2Xk}
(p^2, p v^{(k-1)g_1}, v^{(kp+2)g_1} ).
\end{equation}
This allows us to obtain the description (up to extensions) of the $ku$-homology groups of
$\mz_{p^2} \times \mz_{p^{k+1}}$.

We give constructive proofs, since we are interested in to adapt the results obtained in this work to
other bordism theories, such as, $BP$-homology.

\section{Preliminaries.}\label{s2}

This work is a natural extension of \cite{kuosa}, we will use the same notation.

Let $ku$ be the connective complex $K$-theory. Let $v$ be a generator of $\pi_2 (ku) \simeq \mz.$
The canonical (non typical) Formal Group Law $\cal F$ is given by:
\begin{equation}\label{FGL}
x +_{\cal F}  y = x + y - vxy.
\end{equation}

For a fixed prime number $p$ we will work with the local version of $ku$ and the respective local version of (\ref{FGL}).

For any natural number $n$ we have the formal power series (that is in fact a polynomial):
$$
[p^n](x) = \sum \limits_{k=0}^{p^n-1} a_{k,n} x^{k+1}
$$
where $a_{k,n} = {p^n \choose k+1} v^k$. We omit the dependence on $n$ when it is clear from the context.

For $1 \leq i \leq n$ the term $a_{p^i - 1}$ will be denoted by $a_{g_i} = u_i p^{n-i} v^{g_i}$. Here $u_i$ is a unit in $ku_*$. \\

From here we will use the reduced version of $ku_*$ without further comments. We will denote by
$ku_*(\mz_{p^t} \wedge \mz_{p^n})$ the reduced $ku$-homology of the group $\mz_{p^t} \times \mz_{p^n}$. \\

The element $e_{1,n} \in ku_*(\mz_{p^n})$ (the bottom class) is the so called toral class. In the group
$ku_*(\mz_{p^n} \wedge \mz_{p^t})$ we also have a toral class $\tau$, that comes from the
canonical map $\mz^2 \rightarrow \mz_{p^n} \times \mz_{p^t}$.

We have the K\"unneth map
$\kappa \colon ku_*(\mz_{p^n}) \otimes_{ku_*} ku_*(\mz_{p^t}) \longrightarrow ku_* (\mz_{p^n} \wedge \mz_{p^t}).$
The image of the product of the toral classes $e_{1,n} \otimes e_{1,t}$ under this map is $\tau$. This
map is injective, therefore we have:
$ ann_{ku_*} (e_{1,n} \otimes e_{1,t}) = ann_{ku_*} (\tau).$

We have the Landweber split short exact sequence for bordism theories:
$$
0 \rightarrow ku_*(\mz_{p^t}) \otimes_{{}_{ku_*}} ku_*(\mz_{p^n})  \rightarrow ku_*(\mz_{p^t} \wedge \mz_{p^n}) \rightarrow
\sum {}{\rm Tor}_1^{ku_*} (ku_*(\mz_{p^t}),ku_*(\mz_{p^n})) \rightarrow 0,
$$
therefore we have a direct sum decomposition of the group $ ku_*(\mz_{p^t} \wedge \mz_{p^n})$.
As we did in \cite{kuosa} we approximate the first and the third term of this short exact sequence
by an spectral sequence.

\section{The spectral sequence.}\label{s3}

Let $F$ the free $ku_*$-module in generators $\alpha_i$'s for $i \geq 1$, that is, $F = \bigoplus_{i \geq 1} ku_* \alpha_i.$
For natural numbers $n \geq 2$ and $k \geq n -1$ we consider the map
$$
\begin{array}{rcl}
\partial_{k+1} \colon  F \otimes_{ku_*} ku_*(\mz_{p^n})& \longrightarrow & F \otimes_{ku_*} ku_*(\mz_{p^n}) \\
\alpha_i \otimes e_j & \longrightarrow & \sum \limits_{t=0}^{p^{k+1} -1} a_t \alpha_{i-t} \otimes e_j
\end{array}
$$
where $\alpha_h = 0$ if $h \leq 0$. Here $a_t \in ku_*$ are the coefficients of the $[p^{k+1}]$-series.
The elements $\alpha_i \otimes e_j$ are the basic generators of $ F \otimes_{ku_*} ku_*(\mz_{p^n})$.
Note that
${\rm coker} (\partial_{k+1}) = ku_*(\mz_{p^{k+1}}) \otimes_{ku_*} ku_*(\mz_{p^{n}}).$

We consider the chain complex:
\begin{equation}\label{complex}
\cdots \longrightarrow 0 \longrightarrow F \otimes_{ku_*} ku_*(\mz_{p^n}) \longrightarrow F \otimes_{ku_*} ku_*(\mz_{p^n}) \longrightarrow 0 \longrightarrow \cdots
\end{equation}
where the only non trivial map is given by $\partial_{k+1}$.

Note that every element of the module $F \otimes_{ku_*} ku_*(\mz_{p^n})$
has a unique expression modulo $p^n$, therefore we can define a
filtration:
$$
\begin{array}{ccccc}
|c| = 0 \,\, {\rm for} \,\, c \in ku_*.  & &
|\alpha_i| = i ( p^{k+1} +1 ). & &
|e_j| = j (p^{k} + 1).
\end{array}
$$
We will denote by $c [i,j] = c \, \alpha_i \otimes e_j$ for any $c \in ku_*$.

We will prove that when $n = 2$ and $k\geq2$, in the spectral sequence associated to {\rm (\ref{complex})}, there exist only two families of differentials given (up to units) by:
$$
\begin{array}{rcl}
[i,j] & \longrightarrow & p v^{kg_1} [i, j - k g_1] \\
p[i,j] & \longrightarrow & v^{h(k)} [i-g_1, j - (kp+2) g_1]
\end{array}
$$
Here $h(k) = (k p +2) g_1$. Note that the map we are using is $\partial_{k+1}$.

When $n=e=k+1$ we obtain that there exist $e$ families of diferentials given (up to units) by:
$$
p^t[i,j] \longrightarrow  p^{e-t-1} v^{h(t)} [i-g_t, j - g_{t+1}].
$$
Here $0 \leq t \leq e-1$ and $h(t) = g_t + g_{t+1}$.

\section{Relations in $F \otimes_{ku_*} ku_*(\mz_{p^2})$.}\label{s4}

In this section we construct relations in $F \otimes_{ku_*} ku_*(\mz_{p^2})$ that will be useful to control the
combinatorics in the non filtered differentials. The main result of the section is Theorem~\ref{Teoosa}, this gives an appropiate
expression for $p^k [a,b]$ as an element of $F \otimes_{ku_*} ku_*(\mz_{p^2})$.
In this section the elements $a_i \in ku_*$ are the coeficients of the $[p^2]$-series. We will denote by $a_{g_1} = u_1 p v_{g_1}$.
The sums ${\bf \sum_{1,k+1}^{[a,b]},  \sum_{2,k+1}^{[a,b]} }$, and ${\bf \sum_{3,k+1}^{[a,b]} }$
that appears in the Theorem, in general will be treated inductively of by Smith morphisms arguments.

Recall that for any ordered pair of natural numbers $(i,j)$ with $i,j \geq 0$ we have the Smith morphism $\phi_{i,j}$
given by:
$$
\xymatrix
{
F \otimes_{ku_*} ku_*(\mz_{p^n}) \ar[rr]^{\phi_{i,j}} & & F \otimes_{ku_*} ku_*(\mz_{p^n})  \\
[a,b] \ar[rr] & & [a-i,b-j].
}
$$
These morphisms are compatible with $\partial_t$.

\begin{Lema}\label{1}
For $k \geq 0$ the elements $p^{n+k} [a,b]$ with $1 \leq b \leq (k+1)g_1$ are zero, in the tensor product
$F \otimes_{ku_*} ku_*(\mz_{p^n})$, .
\end{Lema}

\dem We proceed by induction on $k$. The case $k = 0$ is Lemma 5.1 of \cite{kuosa}. Now suppose the result is valid for $0 \leq k$.
For $k+1$ we proceed by induction on b. If $1 \leq b \leq (k+1)g_1$ the result follows from the inductive hypothesis. Now we suppose that
$(k+1)g_1 + 1 \leq b \leq (k+2) g_1$. We assert that
$$
p^{n+k+1}[a,b] = - p^{k+1}\left (  \sum \limits_{i=1}^{p-2}a_{i,n} [a,b-i] \right ).
$$
Here the elements $a_{i,n}$ are the coefficients of the $[p^n]$-series, we will omit the dependence of
$n$ since it is clear from the context. To prove the assertion it suffices to prove that
for $1 \leq j \leq n$ the element $p^{k+1} a_{g_j} [a,b-g_j] = 0$. The case $j=1$ is just the inductive hypothesis; suppose
$j > 1$. Since $g_j < (k+2) g_1$ we have
$$
j < \sum \limits_{i=1}^{j-1} p^i < k+1,
$$
therefore
$$
p^{k+1} a_{g_j} [a, b- g_j] = u p^{n+k_j} v^{g_j} [a, b-g_j].
$$
Here $k_j = (k+1) - j$. Since
$$
b - g_j \leq (k+2)g_1 - g_j
        = g_1 \left ( k+2 - \sum \limits_{i=0}^{j-1}p^i \right )
         \leq g_1 (k+1 - j)
$$
the assertion follows by induction. To deal with the sum
$$
- p^{k+1} \left (   \sum \limits_{i=1}^{p-2} a_i [a, b-i] \right )
$$
we can use an easy inductive argument on $b-i$ since $a_i = u p^n v^i$
where $u \in ku_*$ is a unit.
\QED

For an element $[a,b] \in F \otimes_{ku_*} ku_*(\mz_{p^2})$ let denote by
$$
A^{[a,b]} = \sum \limits_{i=1}^{p-2} {p^2 \choose i+1} v^i [a,b-i] \,\,\,\,\,\,\,\,\,\,\,\,\,\,
B^{[a,b]} = \sum \limits_{i=p}^{p^2-2} {p^2 \choose i+1} v^i [a,b-i]
$$
We also denote by $(A+B)^{[a,b]} = A^{[a,b]} + B^{[a,b]}$.

For a fixed natural number $k \geq 2$ we define $c_k(0) = d_k(0) = 1$,
and for $t \geq 1$ we define
$$
\begin{array}{cccc}
c_k(t) = d_{k-1}(t-1) + c_{k-1}(t) & & & d_k(t) = c_{k-1}(t).
\end{array}
$$
We set $c_2(t) = d_2(t) = 0$ for $t \geq 1$ and $d_k(t) = c_k(t) =
0$ if $t < 0$ or $k < 2$. We also set $v^t = 0$ for $t < 0$. As
usual, for a rational number $q$ we denote by $[q]$ the maximum
integer function. For natural numbers $t,r$ we denote by
$$
\begin{array}{ccccc}
c_t = (-1)^{t+1} u_1^t &  & c_t ' = (-1)^{t+1}c_{2t+1}(t) &  & c_{r,t} = (-1)^{t+r+1} u_1^t c_{t+2r+1}(r).
\end{array}
$$

\begin{Teo}\label{Teoosa}
For an element $[a,b] \in F \otimes_{ku_*} ku_*(\mz_{p^2})$ the following equality holds for $k \geq 2$.
$$
\begin{array}{rl}
p^{k+1}[a,b] & =  (-1)^k u_1^k p v^{kg_1}[a,b-kg_1] \\
\\
            &+  (-1)^k u_1^{k-1}  v^{(k-1)g_1 + g_2}[a,b-(k-1)g_1-g_2] \\
\\
            & + \sum \limits_{t=1}^{\left[\frac{k}{2} \right ]} c_{k+1}(t)(-1)^{k+t} u_1^{k-2t} p v^{(k-2t)g_1 + tg_2}[a,b-(k-2t)g_1 - tg_2] \\
\\
            & + d_{k+1}(t)(-1)^{k+t} u_1^{k-2t-1}  v^{(k-2t-1)g_1 } v^{(t+1)g_2}[a,b-(k-2t-1)g_1 - (t+1) g_2] \\
\\
            & + {\bf \sum_{1,k+1}^{[a,b]} + \sum_{2,k+1}^{[a,b]} +\sum_{3,k+1}^{[a,b]} }
\end{array}
$$
The terms of the last row are given by:
$$
\begin{array}{l}
{\bf \sum_{1,k+1}^{[a,b]}}= \sum \limits_{t=0}^{k-1} c_t p^{k-t-1} v^{tg_1} \left ( A + B \right )^{[a,b-tg_1]},  \\
\\
{\bf \sum_{2,k+1}^{[a,b]}}= \sum \limits_{t=1}^{\left [\frac{k-1}{2} \right ]} c_t' p^{k-2t-1} v^{tg_2} \left ( A + B \right )^{[a,b-tg_2]}, \\
\\
{\bf \sum_{3,k+1}^{[a,b]}}= \sum \limits_{r=1}^{\left [\frac{k-2}{2} \right ]} \sum \limits_{t=1}^{n_{r,k}} c_{r,t} p^{k-2r-t-1} v^{rg_2+tg_1}
\left ( A + B \right )^{[a,b-rg_2-tg_1]}.
\end{array}
$$
Here the limit $n_{r,k} = k-1-2r$.
\end{Teo}

\dem We proceed by induction. For $k=2$ we have
$$
\begin{array}{rcl}
 p^3 [a,b] & = &- p \left ( A^{[a,b]}  + u_1 p v^{g_1} [a, b-g_1] + B^{[a,b]} + v^{g_2} [a,b-g_2] \right )  \\
\\
           & = & - p \left ( A + B \right )^{[a,b]}  + u_1 v^{g_1} \left ( A + B \right )^{[a,b-g_1]} \\
\\
           &   & + u_1^2 p v^{2g_1} [a,b-2g_1] + u_1 v^{g_1+g_2} [a,b-(g_1+g_2)] - p v^{g_2}[a,b-g_2].
\end{array}
$$
Since $c_2(1)=0$ and $v^{t} = 0$ for $t<0$ the result follows for $k=2$. Suppose the result is valid for $k \geq 2$,
we have for $k+2$:
$$
\begin{array}{rcl}
p^{k+2} [a,b] & = & (-1)^k u_1^k p^2 v^{kg_1}[a,b-kg_1] +  (-1)^k u_1^{k-1} p v^{(k-1)g_1 + g_2}[a,b-(k-1)g_1-g_2] \\
\\
            & + & \sum \limits_{t=1}^{\left [\frac{k}{2} \right]} c_{k+1}(t)(-1)^{k-t} u_1^{k-2t} p^2 v^{(k-2t)g_1 + tg_2}[a,b-(k-2t)g_1 - tg_2] \\
\\
            & +  & d_{k+1}(t)(-1)^{k-t} u_1^{k-2t-1} p v^{(k-2t-1)g_1 } v^{(t+1)g_2}[a,b-(k-2t-1)g_1 - (t+1) g_2]  \\
\\
&+&\sum \limits_{t=0}^{k-1} c_t p^{k-t} v^{tg_1} \left ( A + B \right )^{[a,b-tg_1]}  +
\sum \limits_{t=1}^{\left [\frac{k-1}{2} \right ]} c_t' p^{k-2t} v^{tg_2} \left ( A + B \right )^{[a,b-tg_2]}  \\
\\
&+&\sum \limits_{r=1}^{\left [\frac{k-2}{2} \right ]} \sum \limits_{t=1}^{n_{r,k}} c_{r,t} p^{k-2r-t} v^{rg_2+tg_1}
\left ( A  + B \right )^{[a,b-rg_2-tg_1]}.
\end{array}
$$
By the relation imposed in the second factor we have that the expression:
$$
(-1)^k u_1^k p^2 v^{kg_1}[a,b-kg_1]
+ \sum \limits_{t=1}^{\left [\frac{k}{2} \right ]}
c_{k+1}(t)(-1)^{k-t} u_1^{k-2t} p^2 v^{(k-2t)g_1 + tg_2}[a,b-(k-2t)g_1 - tg_2]
$$
is equal to:
$$
\begin{array}{rl}
 & (-1)^{k+1} u_1^{k+1} p v^{(k+1)g_1} [a, b- (k+1)g_1] + (-1)^{k+1} u_1^{k} v^{kg_1+g_2} [a, b-kg_1-g_2] \\
\\
+ & (-1)^{k+1} u_1^k v^{kg_1} \left ( A + B \right )^{[a,b-kg_1]}  \\
\\
+ & \sum \limits_{t=1}^{\left [\frac{k}{2} \right ]} c_{k+1}(t) (-1)^{k+1-t} u_1^{k-2t+1} p v^{(k-2t+1)g_1 +tg_2} [a,b-(k-2t+1)g_1 -tg_2]\\
\\
+ &  c_{k+1}(t) (-1)^{k+1-t} u_1^{k-2t} v^{(k-2t)g_1 +(t+1)g_2} [a,b-(k-2t)g_1 -(t+1)g_2]\\
\\
+ &  c_{k+1}(t) (-1)^{k+1-t} u_1^{k-2t} v^{(k-2t)g_1 +tg_2} \left ( A + B \right )^{[a,b-(k-2t)g_1 - tg_2]}.
\end{array}
$$
Therefore $p^{k+2} [a,b]$ is equal to:
$$
\begin{array}{rl}
  & (-1)^{k+1} u_1^{k+1} p v^{(k+1)g_1} [a, b- (k+1)g_1] + (-1)^{k+1} u_1^{k} v^{kg_1+g_2} [a, b-kg_1-g_2] \\
\\
+ & \sum \limits_{t=1}^{\left [\frac{k}{2} \right ]} \left ( c_{k+1}(t)+d_k(t-1) \right )
(-1)^{k+1-t} u_1^{k-2t+1} p v^{(k-2t+1)g_1 +tg_2} [a,b-(k-2t+1)g_1 -tg_2] \\
\\
+ & c_{k+1}(t) (-1)^{k+1-t} u_1^{k-2t} v^{(k-2t)g_1 +(t+1)g_2} [a,b-(k-2t+1)g_1 -tg_2]\\
\\
+ & (-1)^{k+1} u_1^k v^{kg_1} \left ( A + B \right )^{[a,b-kg_1]} \\
\\
+ & \sum \limits_{t=1}^{\left [\frac{k}{2} \right ]}
c_{k+1}(t) (-1)^{k+1-t} u_1^{k-2t} v^{(k-2t)g_1 +tg_2} \left ( A + B \right )^{[a,b-(k-2t)g_1 - tg_2]} \\
\\
+ &p {\bf \sum_{1,k+1}^{[a,b]} } + p {\bf \sum_{2,k+1}^{[a,b]} } + p {\bf \sum_{3,k+1}^{[a,b]} } \\
\\
+ & (-1)^{k - \left [\frac{k}{2} \right ]} d_{k+1} \left( \left [\frac{k}{2} \right ] \right ) u_1^{k-2\left [\frac{k}{2} \right] -1}
p v^{(k-2[\frac{k}{2}]-1) g_1} v^{\left( \left [\frac{k}{2} \right ]+1 \right ) g_2} [a,b-(k-2t-1)g_1 - (t+1)g_2 ]
\end{array}
$$
Note that the last term in the previous expression
is not zero only when $k$ is odd; this is the reason because we separated the powers of $v$ in the statement
of the theorem. This last term, in general, takes care of the cases that arise from the parity of $k$.

If we take the first $\left [\frac{k}{2} \right ] -1$ terms of the sum
\begin{equation}\label{suma4}
\sum \limits_{t=1}^{\left [\frac{k}{2} \right ]}
c_{k+1}(t) (-1)^{k+1-t} u_1^{k-2t} v^{(k-2t)g_1 +tg_2} \left ( A + B \right )^{[a,b-(k-2t)g_1 - tg_2]}
\end{equation}
and we add with $p {\bf \sum_{3,k+1}^{[a,b]}}$, we obtain
$$
\sum \limits_{r=1}^{\left [\frac{k-2}{2} \right ]} \sum \limits_{t=1}^{n_{r,k+1}} c_{r,t} p^{k-2r-t} v^{rg_2+tg_1}
\left ( A  + B \right )^{[a,b-rg_2-tg_1]}.
$$
Now, if $k$ is even the last term of (\ref{suma4}) is given by:
$$
(-1)^{\frac{k}{2} + 1} c_{k+1} \left ( \frac{k}{2} \right ) v^{\frac{k}{2} g_2} \left ( A + B \right )^{[a, b-\frac{k}{2} g_2]}.
$$
Adding this term with $p {\bf \sum_{2,k+1}^{[a,b]}}$ we obtain  ${\bf \sum_{2,k+2}^{[a,b]}}$.
When $k$ is odd, the last term of (\ref{suma4}) is equal to:
$$
(-1)^{\frac{k-1}{2}} u_1 c_{k+1} \left ( \frac{k-1}{2} \right )
v^{g_1+\frac{k-1}{2} g_2} \left ( A + B \right )^{[a, b-g_1-\frac{k-1}{2} g_2]}.
$$
Since $n_{\frac{k-1}{2},k+1} = 1$ we can add this term to $p {\bf \sum_{3,k+1}^{[a,b]}}$, and the theorem follows
from the fact that $c_{\frac{k-1}{2},1} = (-1)^{\frac{k-1}{2}} u_1 c_{k+1} \left (\frac{k-1}{2} \right )$.
\QED

\begin{Remark}\label{ck} \rm From the inductive definition of the functions
$c_k(t)$ and $d_k(t)$ and from the fact that $c_2(t) = d_2(t) = 0$
for $t\geq 1$ it follows that for $k$ even:
$$
\begin{array}{l}
c_k(t)= 0 \,\,\,\,\,\, {\rm and} \,\,\,\,\,\, d_k(t)=0 \,\,\,\,\,\,
{\rm for} \,\,\,\,\,\, t
\geq \frac{k}{2}.   \\
\end{array}
$$
and
$$
\begin{array}{l}
 c_{k+1} \left( \frac{k}{2}  \right) =
d_{k}\left(  \frac{k}{2} -1 \right) \\
\\
d_{k+1} (\frac{k}{2}) = 0 \\
\\
c_{ k+1} (t) = 0  \,\,\,\,\,\, {\rm for}  \,\,\,\,\,\, t \geq \frac{k}{2} + 1  \\
\\
d_{k+1} (t) = 0 \,\,\,\,\,\, {\rm for}  \,\,\,\,\,\, t \geq
\frac{k}{2} + 1
\end{array}
$$
\end{Remark}

\begin{Remark}\rm
It is possible to have a compact presentation of these results but we
are interested in to make accessible the material for the non
experts. On the other hand, using this presentation it is possible
to adapt these calculations to other homology theories. In fact we
believe that we can make slight improvements of the formulas
obtained in this work to collapse the respective spectral sequence
for $BP$.
\end{Remark}

Now we are interested in construct a relation in $F \otimes_{ku_*} ku_*(\mz_{p^2})$ in order to give potential differentials
in the spectral sequence associated to {\rm (\ref{complex})} for $k \geq 2$. We begin defining
inductively a family of polynomials with coefficients in $\mz_{(p)}$. We define
$$
\begin{array}{rl}
p_0^k (t) & = d_k(t) + d_{k-1} (t) \\
           & = a_{0,1} d_k(t) + a_{0,2} d_{k-1} (t). \\
           \\
p_1^k (t) & = p_0^k (0) (d_k(t) + d_{k-1}(t)) + (p_0^k (0) - a_{0,1} ) d_{k-2}(t) \\
          & = a_{1,1} d_k(t) + a_{1,2} d_{k-1} (t) + a_{1,3} d_{k-2}(t).
\end{array}
$$
In general we will denote by $a_{i,l}$ the coefficient of $d_{k-l+1} (t)$ in the polynomial $p_i^k(t)$:
$$
p_i^k (t) = a_{i,1} d_k(t) + a_{i,2} d_{k-1} (t) + \cdots + a_{i,i+2} d_{k-i-1} (t).
$$
The coefficients $a_{i,l}$ are given by:
$$
\begin{array}{ccc}
 a_{i,1} = p_{i-1}^k (0) & & a_{i,2} = p_{i-1}^k (0),
\end{array}
$$
and for $3 \leq l \leq i+2$ we define $a_{i,l} = a_{i,l-1} - a_{i-1,l-2}$.

For fixed $k$ and $[a,b] \in F \otimes_{ku_*} ku_*( \mz_{p^2})$ we
denote by ${\cal S}_k^{[a,b]}$ the term:
$$
p^{k+2} [a,b] \,\, + \,\, p^{k+1} \left ( u_1 v^{g_1} [a,b-g_1] -
u_1^{-1} v^{pg_1} [a,b-pg_1] \,\, - \,\, u_1^{-3}  v^{f_{-1} g_1}
[a,b-f_{-1} g_1] \right ).
$$

Here $f_i = (i+3) g_1 + 1$ for $i \geq -3$.

It is not difficult to verify using Theorem~\ref{Teoosa} that as an
element of the tensor product $F \otimes_{ku_*} ku_*(\mz_{p^2})$ the
term ${\cal S}_k^{[a,b]}$ is equal to:
$$
\begin{array}{l}
 \sum \limits_{t=0}^{\left[ \frac{k-2}{2} \right ] }
d_{k,k-1}(t) (-1)^{k-t} u_1^{k-2t-5} p v^{(k+tg_1+f_0)g_1} [a,b-(k+tg_1+f_0)g_1] \\
\\
+ d_{k-1,k-2}(t) (-1)^{k+t} u_1^{k-2t-6}  v^{(k-2t-3)g_1}
v^{(t+4)g_2 -3g_1} [a,b-(k+(t+4)g_1 + 2) g_1] + {\Large \bf \sum}
\end{array}
$$
Here $d_{k,k-1} (t) = d_k(t) + d_{k-1} (t)$ and ${\Large \bf \sum}$
is the sum of the terms ${\bf \sum_{i,k+2}^{[a,b]}}$ for $i = 1,2,3$
and the respective terms that arise from the multiples of $p^{k+1}$.
In general we will ignore this terms since they can be treated by
Smith morphism arguments.

For $0 \leq n \leq k-3$ we define the element in $F \otimes_{ku_*}
ku_*(\mz_{p^2})$:
$$
{\cal S}^{[a,b]}_{k,n} = \sum \limits_{i=0}^{n} p_i^k (0)
u_1^{-(2i+5)} v^{f_i g_1} [a,b-f_i g_1]
$$

\begin{Propo}\label{Sdif}\rm
In the tensor product $F \otimes_{ku_*} ku_*(\mz_{p^2})$ for $0 \leq
n \leq k-5$ the element ${\cal S}^{[a,b]}_{k} - p^{k+1} {\cal
S}^{[a,b]}_{k,n}$ is equal to:
$$
\begin{array}{l}
\sum \limits_{t=0}^{\left [ \frac{k}{2} -1 \right ]}  (-1)^{k+t} p_{n+1}^k (t)
u_1^{k-(2n+2t+7)} p v^{f(k,t,n)}  [a,b-f(k,t,n)] \\
\\
 +   (-1)^{k+t}  p_{n+1}^{k-1} (t) u_1^{k-(2n+2t+8)} v^{h_1(k,t)} v^{h_2(t,n)} [a,b-h_{1,2}(k,n,t)] + {\Large \bf \sum}.
\end{array}
$$
Here $f(k,t,n) = (k + (t+n+4)g_1+1) g_1$, $h_1(k,t) = (k-2t-3)g_1,
\,\, h_2(t,n)= (t+2)g_2 + f_n g_1$ and $h_{2,1}(k,n,t) = h_1(k,t) +
h_2 (t,n)$. The element $\Large \bf \sum$ is the sum of the terms
${\bf \sum_{i,k+2}^{[a,b]}}$ for $i = 1,2,3$ and the respective
terms that arise from the multiples of $p^{k+1}$.
\end{Propo}

\dem We proceed by induction on $n$. For $n=0$ we have that ${\cal S} - p^{k+1} {\cal S}_0$ is equal to:
$$
\begin{array}{l}
\sum \limits_{t=0}^{\left[ \frac{k-2}{2} \right ]}
d_{k,k-1}(t) (-1)^{k+t} u_1^{k-2t-5} p v^{(k+tg_1 + f_0) g_1} [a, b - (k + tg_1 + f_0) g_1] \\
\\
+ d_{k-1,k-2}(t) (-1)^{k+t} u_1^{k-2t-6}  v^{(k-2t-3)g_1}
v^{(t+4)g_2 - 3 g_1}
[a,b - (k + tg_1 + f_0 + p) g_1] + {\Large \bf \sum} \\
\\
+ \sum \limits_{t=0}^{\left [ \frac{k}{2} \right ]}
p_0^k (0) c_{k+1}(t) (-1)^{k+1+t} u_1^{k-2t-5} p v^{(k-2t-5)g_1 + (t+3)g_2} [a,b - ( k+1 + (t+3) g_1 )g_1] \\
\\
+ p_0^k (0) d_{k+1}(t) (-1)^{k+1+t} u_1^{k-2t-6} v^{(k-2t-1)}g_1
v^{(t+4)g_2 -5g_1} [a,b - ( k+1 + (t+3) g_1 + p )g_1 ]
\end{array}
$$
The terms of the sums when $t=0$ are canceled. Since we have that
$d_k(t) + d_{k-1}(t) - p_0^k c_{k+1}(t)$ is equal to:
$$
\begin{array}{rl}
& c_{k-1} (t) + c_{k-2} (t) - p_0^k(0) (d_k(t-1) + c_{k-1}(t) + d_{k-1}(t-1)) \\
\\
= &c_{k-2} (t)  - p_0^k(0) (d_k(t-1) + d_{k-1}(t-1)) - c_{k-1}(t) \\
\\
                                      = & - p_0^k(0) (d_k(t-1) + d_{k-1}(t-1)) - d_{k-2} (t-1)\\
\\
 = & -p_1^k(t-1).
\end{array}
$$
This gives the terms for $0 \leq t \leq \left [ \frac{k}{2} -2
\right ]$.

Now let see the behavior when $t = \left [ \frac{k}{2} -1 \right ]$.
We will only consider the case when $k$ is even; the case when $k$
is odd is proved using exactly the same argument. We have by
Remark~\ref{ck} the following equality:
$$
p_0^k(0) c_{k+1} \left ( \frac{k}{2} \right ) = p_0^k(0) d_k \left (
\frac{k}{2}  - 1 \right ).
$$
We have by definition that:
$$
p_1^k \left (  \frac{k}{2}  - 1\right ) = p_0^k(0) \left (  d_k
\left ( \frac{k}{2}  - 1 \right  ) + d_{k-1} \left (  \frac{k}{2} -1
\right ) \right ) + (p_0^k(0) - a_{0,1}) d_{k-2} \left ( \frac{k}{2}
- 1 \right ).
$$
But $d_{k-2} \left (  \frac{k}{2} - 1 \right ) =  d_{k-1} \left (
\frac{k}{2} - 1 \right ) = 0 $. On the other hand we have:
\[
\begin{array}{rl}
p_0^k(0) d_{k+1} \left ( \frac{k}{2} - 1 \right ) & = p_0^k(0) c_k
\left ( \frac{k}{2}  - 1 \right ) \\
\\
& =  p_0^k(0) \left ( c_{k-1} \left ( \frac{k}{2}  - 1 \right )  +
 d_{k-1} \left ( \frac{k}{2}  - 2 \right ) \right )\\
 \\
& = p_0^k(0) \left ( d_{k-1} \left ( \frac{k}{2}  - 2 \right )  +
 d_{k-2} \left ( \frac{k}{2}  - 2 \right ) \right ).
\end{array}
\]
We also have:
$$
p_1^{k-1} \left (  \frac{k}{2}  - 2  \right ) = p_0^{k-1}(0) \left (
d_{k-1} \left ( \frac{k}{2}  - 2 \right  ) + d_{k-2} \left (
\frac{k}{2} - 2 \right ) \right ) + (p_0^{k-1} (0)- a_{0,1}) d_{k-3}
\left ( \frac{k}{2} - 2 \right ).
$$
Since $d_{k-3} \left (  \frac{k}{2} - 2 \right ) = 0$, the result
follows for $n=0$. Now suppose the result is valid for $0 \leq n
\leq k-6$; for $n+1$ we have that the term:
$$
- p^{k+1} p_{n+1}^k (0) u_1^{-(2n+7)} v^{f_{n+1} g_1} [a,b-f_{n+1}
g_1]
$$
is equal to:
$$
\begin{array}{l}
\sum \limits_{t=0}^{\left [ \frac{k}{2} \right ]} p_{n+1}^k (0) c_{k+1}(t) (-1)^{k+1+t}
u_1^{k-2t-2n-7} p v^{(k-2t+f_{n+1})g_1 +tg_2} [a, b - (k-2t+f_{n+1}) g_1 - t g_2 ] \\
\\
+ p_{n+1}^k (0) d_{k+1}(t) (-1)^{k+1+t} u_1^{k-2t-2n-8}
v^{(k-2t-1)g_1} v^{F_{n+1}(t)} [a,b-(k-2t-1) g_1 - F_{n+1}(t)]
\end{array}
$$
Here $F_{n+1} (t)  = f_{n+1} g_1 + (t+1) g_2$. We have that
$p_{n+1}^k (0) c_{k+1}(t) - p_{n+1}^k (t)$ is by definition:
$$
\left ( \sum \limits_{i=1}^{n+3} a_{n+1,i} \right ) c_{k+1}(t) - \sum \limits_{i=1}^{n+3} a_{n+1,i} d_{k+1-i} (t).
$$
We also have that $a_{n+1,1} (c_{k+1}(t) - d_k(t)) = a_{n+1,1}
(d_k(t-1) + d_{k-1} (t-1) )$, and for $1 \leq i \leq n+2$:
$$
a_{n+1,i+1} (c_{k+1} (t) - d_{k-i} (t)) = a_{n+1,i+1} (d_k(t-1) + d_{k-1} (t-1) + \cdots + d_{k-i-1} (t-1)).
$$
Therefore we obtain that  $p_{n+1}^k (0) c_{k+1}(t) - p_{n+1}^k (t)
= p_{n+2}^k (t-1)$. This takes care of the terms with $0 \leq t \leq
\left [ \frac{k}{2} -2 \right ]$. Finally note that for $k$ even:
$$
\begin{array}{rl}
p_{n+1}^k (0) d_{k+1} (\frac{k}{2} - 1) & = p_{n+1}^k (0) \left (
c_{k-1} (\frac{k}{2} - 1) + d_{k-1} (\frac{k}{2} - 2 ) \right ) \\
\\
& = p_{n+1}^k (0) \left ( d_{k-1} (\frac{k}{2} - 2) + d_{k-2}
(\frac{k}{2} - 2 ) \right ),
\end{array}
$$
since $c_{k-2} (\frac{k}{2} -1) = 0$. We also have that:
$$
p_{n+2}^{k-1} \left (  \frac{k}{2} - 2 \right ) = p_{n+1}^{k-1} (0)
\left (  d_{k-1} \left ( \frac{k}{2} - 2\right ) + d_{k-2} \left (
\frac{k}{2} - 2 \right) \right ),
$$
since $d_{k-3} \left (  \frac{k}{2} - 2 \right ) = 0$. We have by
hypothesis $n \leq 5$, this implies that $p_{n+1}^k (0) =
p_{n+1}^{k-1} (0)$. On the other hand
$$
\begin{array}{rl}
p_{n+1}^k (0) c_{k+1} \left (  \frac{k}{2} \right ) & = p_{n+1}^k(0)
d_k \left ( \frac{k}{2} - 1  \right ) \\
\\
& = p_{n+2}^k (\frac{k}{2} - 1)
\end{array}
$$
since $c_k \left ( \frac{k}{2} \right ) = 0$ and $ d_{k-1} \left (
\frac{k}{2} - 1 \right ) = 0$. This completes the proof for $k$
even.
\QED

\section{Calculating the differentials for $\mz_{p^k} \times \mz_{p^2}$.}\label{s5}

In this section we construct in Theorem~\ref{dife} the two potential differentials of the spectral sequence of
the group $\mz_{p^k} \times \mz_{p^2}$. In Proposition~\ref{CP} we prove that the elements that can give rise
to other differentials are in fact permanent cycles. This will prove that (\ref{an2Xk}) is the annihilator
of the $ku$-toral class for this group. Also it will give us a presentation (up to extensions)
for the $ku$-homology of the group $\mz_{p^k} \times \mz_{p^2}$.

We will denote the elements of the $[p^2]$-series by $a_i = w_i p^2 v^i$ for $0 \leq i \leq p^2-2$ with $i \neq g_1$
and by $a_{g_1} = u_1 p v^{g_1}$. Here $u_1$ is a unit in $ku_*$.
Recall that $a_{g_2} = v^{g_2}$. For $k \geq 2$ we will denote the elements of the $[p^{k+1}]$-series
by $a_{i,k+1} = z_i p^{k+1-t} v^i$ if $g_t < i < g_{t+1}$ for $1 \leq t \leq k+1$; we will denote $a_{g_t,k+1} = y_t p^{k-t} v^{g_t}$.
Here $y_t$ is also a unit of $ku_*$. Recall that for any natural number $k$ the element $a_{0,k}$ of the $[p^k]$-series is given by $p^k$.

\begin{Remark}\label{sumas}
It is not difficult to verify that ${\bf \sum_{2,k+1}^{[a,b]} }$ and
${\bf \sum_{3,k+1}^{[a,b]} }$ are Smith morphism images of ${\bf
\sum_{1,k+1}^{[a,b]}}$. On the other hand we have that:
$$
\begin{array}{rl}
{\bf \sum_{1,k+1}^{[a,b]}}&  = \sum \limits_{t=0}^{k-1}  \sum \limits_{i=1}^{p-2} c_t w_i p^{k-t+1} v^{i+tg_1} [a,b-(tg_1+i)] \\
\\
& +  \sum \limits_{t=0}^{k-1}  \sum \limits_{j=1}^{pg_{{}_1}-1} c_t w_{j+g_1} p^{k-t} v^{j+(t+1)g_1} [a,b-(t+1)g_1-j]
\end{array}
$$
Therefore, in order to prove that ${\bf \sum_{1,k+1}^{[a,b]}} + {\bf \sum_{2,k+1}^{[a,b]}} + {\bf \sum_{3,k+1}^{[a,b]}}$
is an element of ${\rm Im} \partial_{k+1}$ we only consider the element:
$$
 \sum \limits_{t=0}^{k-1} c_t w_1 p^{k-t+1} v^{tg_1 +1} [a,b- (tg_1+1)] + c_{k-1} w_p p v^{kg_1 +1} [a,b-(kg_1 +1)].
$$
\end{Remark}

The following result gives a family of elements that are on ${\rm Im} \, \partial_{k+1}$.

\begin{Lema}\label{nuevo}
The elements $p v^{kg_1} [a,b]$ with $1 \leq a \leq g_2$ and $1 \leq
b \leq g_1^2$ are on ${\rm Im} \, \partial_{k+1}$. Moreover, we have
(up to units):
$$
\partial_{k+1} \left (  [a,b+kg_1]  + \low \right ) = p v^{kg_1} [a,b].
$$
Here, as usual, ``$\low$'' stands lower filtration terms than $[a,b+kg_1]$.
\end{Lema}

\dem The fact that $p v^{kg_1} [1,1] \in \, {\rm Im} \, \partial_{k+1}$ is item $a)$ of Theorem~\ref{dife}. Now suppose we have proved the assertion for
$p v^{kg_1} [1,b]$ with $b < g_1^2$. We have that
$$
\begin{array}{rl}
\partial_{k+1} ([1,b+1+kg_1]) & = p^{k+1} [1,b+1+kg_1] \\
\\
 & = (-1)^k u_1^k p v^{kg_1} [1,b+1] + {\bf \sum_{1,k+1}^{[a,b]}} + {\bf \sum_{2,k+1}^{[a,b]}} + {\bf \sum_{3,k+1}^{[a,b]}}
\end{array}
$$
By the previous remark we only analize:
$$
\sum \limits_{t=0}^{k-1} c_t w_1 p^{k-t+1} v^{tg_1 +1} [1,(k-t)g_1 + b] + c_{k-1} w_p p v^{kg_1 + 1} [a,b].
$$
The last term is in ${\rm Im} \, \partial_{k+1}$ by the inductive hypothesis. For $0 \leq t \leq k-1$ we proceed by
reverse induction on $t$ in order to prove that
\begin{equation}\label{indt}
 p^{k-t+1} v^{tg_1 +1} [1,(k-t)g_1 + b] \,\, \in \,\, {\rm Im} \,\, \partial_{k+1}.
\end{equation}
For $t = k-1$ we have that
$$
\begin{array}{rl}
p^2 v^{(k-1)g_1 + 1} [1, g_1 + b] & = - \sum \limits_{i =0}^{p-2} w_i p^2 v^{(k-1)g_1+i+1} [1,b+g_1-i]  - u_1 v^{kg_1+1} p [1,b]\\
\\
& - \sum \limits_{j=1}^{g_1^2} w_{j+g_1} p v^{kg_1 + j} [1,b-j].
\end{array}
$$
By the inductive hypothesis on the second coordinate, we only have to consider
$$
- \sum \limits_{i =0}^{p-2} w_i p^2 v^{(k-1)g_1+i+1} [1,b+g_1-i]  - u_1 v^{kg_1+1} p [1,b].
$$
Therefore, by a Smith morphism argument the assertion follows for $t = k-1$.
Suppose we have proved that for $0 < t \leq k-1$ the term (\ref{indt}) is an element of ${\rm Im} \, \partial_{k+1}$.
For $t-1$ we have that the element $p^{k-t+2} v^{(t-1)g_1 +1} [1,(k-t+1)g_1 + b]$ is equal to:
$$
 (-1)^{k-t+1} u_1^{k-t+1} p v^{kg_1+1} [1,b] +
{\bf \sum {}_{1,k-t+2}^{[1,b']}} + {\bf \sum {}_{2,k-t+2}^{[1,b']}} + {\bf \sum {}_{3,k-t+2}^{[1,b']}}
$$
Here $b' = (k-t+1)g_1 + b$. Again, we only consider the sum
$$
\sum \limits_{t_1=0}^{k-t} p^{k-t-t_1+2} v^{(t+t_1-1)g_1 + 2} [1.b+(k-t-t_1+1)g_1-1].
$$
We only have to verify that $p^{k-t+2} v^{(t-1)g_1 + 2} [1, b + (k-t+1)g_1 -1] \in  {\rm Im} \, \partial_{k+1}$.
This follows by an inductive argument on the second coordinate.
Suppose we have proved the Lemma for $1 \leq a < g_2$ and $1 \leq b \leq g_1^2$. We have that $\partial_{k+1} ([a+1,a+kg_1])$
is equal to:
$$
\sum \limits_{i=0}^{p-2}  z_i p^{k+1} v^i [a+1-i, 1+kg_1] + y_1 p^k v^{g_1} [a+1-g_1,1+kg_1]
+ \sum \limits_{i=p}^{p^2-2} z_i p^k v^i [a+1-i,1+kg_1].
$$
By a Smith morphism argument we will only prove that:
$$
p^{k+1} [a+1,1+kg_1] \in \, {\rm Im} \, \partial_{k+1}
\,\,\,\, {\rm and} \,\,\,\, p^k v^{g_1} [a+1-g_1, 1+kg_1] \in \, {\rm Im} \, \partial_{k+1}.
$$
We have that $p^{k+1} [a+1,1+kg_1]$ is equal to:
$$
(-1)^k u_1^k p v^{kg_1} [a+1,1] + {\bf \sum {}_{1,k+1}^{[a+1,1+kg_1]}} + {\bf \sum {}_{2,k+1}^{[a+1,1+kg_1]}} + {\bf \sum {}_{3,k+1}^{[a+1,1+kg_1]}}.
$$
A reverse induction argument on $t$ proves that each term of the sum
$$
\sum \limits_{t=0}^{k-t} p^{k-t+1} v^{tg_1 + 1} [a,(k-t)g_1]
$$
is the image under $\partial_{k+1}$ of a ``$\low$'' filtration term. We have that the element  is equal to:
$$
\begin{array}{rl}
p^k v^{g_1} [a+1-g_1, 1+kg_1] & = (-1)^{k-1} u_1^{k-1} p v^{k g_1} [a+1-g_1,1+g_1] \\
\\
&+ c_k(1) (-1)^{k-2} u_1^{k-3} p v^{(k-1)g_1+g_2} [a+1-g_1, 1+3g_1-g_2] \\
\\
& + v^{g_1} \left ( {\bf \sum {}_{1,k}^{[a',b']}} + {\bf \sum {}_{2,k}^{[a',b']}} + {\bf \sum {}_{3,k}^{[a',b']}} \right ).
\end{array}
$$
Here $a' = a+1-g_1$ and $b'=1+kg_1$. By the inductive hypothesis, we have only to verify that
$$
v^{g_1} \sum \limits_{t=0}^{k-2} c_t w_1 p^{k-t} v^{tg_1+1} [a+1-g_1, (k-t)g_1].
$$
It is easy to verify that each summand is on ${\rm Im} \, \partial_{k+1}$ and it is the
image of a ``$\low$'' filtration term. \QED

\begin{Coro}\label{nuevo2}
We have up to units:
$$
\partial_{k+1} ([p, g_1^2 + k g_1 +1] + u_1^{-1} y_1 [1, g_1^2 + (k+1) g_1 +1])
= p v^{k g_1} [p,g_1^2 +1] + \partial_{k+1}(\low).
$$
\end{Coro}

\begin{Remark}\label{bla}\rm
We have by Proposition~\ref{Sdif} that if $k \geq 5$, $a=p$ and $b=(p+1) g_1 +1$
the term ${\mathcal S}_k^{[a,b]} + p^{k+1} {\mathcal S}_{k,k-5}^{[a,b]}$
is equal to:
$$
\begin{array}{rl}
\sum \limits_{t=0}^{\left [ \frac{k}{2} -1 \right]}
& (-1)^{k+t} p_{k-4}^{k} (t) u_1^{-k-2t+3} p v^{f(k,t,k-5)} [p, b - f] \\
& (-1)^{k+t} p_{k-4}^{k-1} (t) u_1^{-k-2t+2}  v^{h_1(k,t)} v^{(h_2(t,k-5))} [p, b - h_{1,2}] + \Sigma.
\end{array}
$$
It is not difficult to verify that this sum is equal to:
$$
(-1)^k p_{k-4}^{k} (0) u_1^{-k+3} p v^{f(k,0,k-5)} [p, 1+g_1^2] + (-1)^{k+1} p_{k-4}^{k} (1) u_1^{-k+1}
p v^{f(k,1,k-5)} [p, 1]  + \Sigma.
$$
This last term is zero modulo $\partial_{k+1}(\low)$ by Lemma~\ref{nuevo}.
\end{Remark}

\begin{Teo}\label{dife}
For $k \geq 2$ we have (up to units) that:
\begin{itemize}
\item[a)] $\partial_{k+1} ([1,kg_1 +1]) = p v^{kg_1} [1,1] $
\item[b)] $\partial_{k+1} (p [p,(pk+1)g_1 + 1] + \low) = v^{(kp+2)g_1} [1,1]$
\end{itemize}
\end{Teo}

\dem To prove item a) we only have to verify that:
$$
{\bf \sum {}_{1,k+1}^{[1,1+kg_1]}} + {\bf \sum
{}_{2,k+1}^{[1,1+kg_1]}} + {\bf \sum {}_{3,k+1}^{[1,1+kg_1]}} = 0,
$$
as an element of $F \otimes_{ku_*} ku(\mz_{p^2})$. By
Remark~\ref{sumas} it suffices to prove that:
$$
 \sum \limits_{t=0}^{k-1} c_t w_1 p^{k-t+1} v^{tg_1 +1} [1,(k-t)g_1] =0.
$$
By Lemma~\ref{1} we have that:
$$
p^{k-t+1} v^{tg_1 + 1} [1,(k-t)g_1] = 0,
$$
and a) follows.

To prove item b) for $k=2$ we first verify that the following expression is zero
modulo $\partial_3(\low)$.
$$
\begin{array}{l}
- p^2 (A + B)^{[p,(2p+1) g_1 +1]} + u_1^{-1} v^{pg_1} p (A+B)^{[p,g_2+1]} \\
\\
+ u_1^{-3} v^{2g_2 - 3 g_1} p (A+B)^{[p,2g_1+1]} - u_1^{-2} v^{2pg_1} (A+B)^{[p,g_1+1]}.
\end{array}
$$
We only have to consider the first two terms since Lemma~\ref{nuevo} takes care
of the last two. For $1 \leq j \leq p-2$ we denote by $z_j p^3 v^j$ the $j$-th coefficient of
the $[p^3]$-series and by $y_1 p^2 v^{g_1}$ the $g_1$-th coefficient of this series.
Note that the elements $z_j$ and $y_1$ are units in $ku_*$. We have that
$$
\begin{array}{rl}
\sum \limits_{i=1}^{p-2} \partial_3 (w_i p v^i [p, (2p+1)g_1 +1-i]) &=  \sum \limits_{i=1}^{p-2} w_i p^4 v^i
[p,(2p+1)g_1 +1 -i] \\
\\
& + \sum \limits_{i=1}^{p-2} \sum \limits_{i_1=1}^{p-2} w_i z_{i_1} p^4 v^{i+i_1} [p-i_1,(2p+1)g_1 + 1 -i] \\
\\
& + \sum \limits_{i=1}^{p-2} w_i y_1 p^3 v^{i+g_1} [1,(2p+1)g_1 +1-i].
\end{array}
$$
and
$$
- \sum \limits_{i=1}^{p-2} \partial_3 (w_1 y_1 v^{i+g_1} [1, (2p+1)g_1+1-i ]) =
-\sum \limits_{i=1}^{p-2} w_i y_1 v^{i+g_1} p^3 [1, (2p+1) g_1 +1].
$$
Therefore modulo $\partial_3(\low)$ we can replace $-p^2 A^{[p,(2p+1) g_1 + 1]}$ by
$$
\sum \limits_{i=1}^{p-2} \sum \limits_{i_1=1}^{p-2} w_i z_{i_1} p^4 v^{i+i_1} [p-i_1,(2p+1) g_1 + 1-i].
$$
Now an easy inductive argument proves that this expression is zero modulo $\partial_3(\low)$.
On the other hand we have that:
$$
\sum \limits_{i=1}^{pg_1-1} \partial_3 (w_{j+g_1} v^{g_1+j} [p, 2pg_1 +1-j])
$$
is equal to
$$
\begin{array}{l}
\,\,\,\,\, \sum \limits_{j=1}^{pg_1 -1} p^3 w_{j+g_1} v^{g_1+j} [p, 2pg_1 +1-j] +
 \sum \limits_{j=1}^{pg_1 -1} \sum \limits_{i_1=1}^{p-2} p^3 w_{j+g_1} z_{i_1} v^{g_1+j+i_1} [p-i_1, 2pg_1 +1-j] \\
 \\
 + \sum \limits_{j=1}^{pg_1 -1} w_{j+g_1} y_1 p^2 v^{2 g_1+j} [1,2pg_1 +1-j].
\end{array}
$$
To deal with the last row of this expression we only consider those terms that are not
divisible by $p^3$. It is not difficult to verify that the sum:
\begin{equation}\label{aux2}
\sum \limits_{h=1}^{p-2} y_1 w_{hp+g_1} p^2 v^{j+2g_1} [1, 2pg_1 +1-j]
\end{equation}
is equal to:
$$
\sum \limits_{h=1}^{p-2} x_h p v^{ph+3g_1} [1, (g_1-h) p +1] +
\sum \limits_{h=1}^{p-2} \sum \limits_{t=1}^{g_{{}_{1}}^2} x_{t,h} p v^{ph+3g_1+t} [1, (g_1-h) p +1 -t].
$$
This last expression is zero modulo $\partial_3(\low)$ by Lemma~\ref{nuevo}. Therefore applying an inductive argument
in the first factor, we have proved that the term $-p^2 B^{[p,(2p+1) g_1 +1]}$ is zero modulo $\partial_3(\low)$.
Now we want to treat the term $u_1^{-1} p v^{pg_1} A^{[p,g_2+1]}$; for to do this we consider the
following:
$$
\sum \limits_{i=1}^{p-2} \partial_3 (v^{pg_1 + i} w_i [p, g_2 +1-i] ).
$$
This expression is equal to:
$$
\begin{array}{l}
\,\,\,\,\, \sum \limits_{i=1}^{p-2} w_i p^3 v^{pg_1 + i} [p, g_2 +1-i]
+ \sum \limits_{i=1}^{p-2}\sum \limits_{i_1=1}^{p-2} w_i z_{i_1} p^3 v^{pg_1 + i + i_1} [p-i_1, g_2 +1-i]
\\
+ \sum \limits_{i=1}^{p-2} w_i y_1 p^2 v^{g_2 + i} [1, g_2 +1-i].
\end{array}
$$
The last row of the previous formula es equal to:
$$
\sum \limits_{i=1}^{p-2} x_i p v^{g_2+g_1+i} [1, pg_1 +1-i]
+ \sum \limits_{i=1}^{p-2} \sum \limits_{t=1}^{pg_1-i} x_{i,t} p v^{g_2+g_1+t+i} [1, pg_1 +1-(i+t)].
$$
All the terms of the sum are Smith morphism images of the first one (we will ignore the coefficients
$x_i$ and $x_{i,t}$). Therefore we only analyze the element
$p v^{g_2+g_1+1} [1,pg_1]$. We have that
$$
\begin{array}{rl}
\partial_3 (v^{pg_1 +1} [1,(p+2) g_1]) & = p^3 v^{pg_1 +1} [1,(p+2) g_1] \\
\\
& = u_1^2 p v^{(p+2)g_1 + 1} [1,pg_1] - c_3(1) p v^{g_2 + pg_1 + 1} [1,g_1] \\
\\
& + \, v^{pg_1 + 1} \left (
{\bf \sum {}_{1,3}^{[1,(p+2)g_1]}} + {\bf \sum
{}_{2,3}^{[1,(p+2)g_1]}} + {\bf \sum {}_{3,3}^{[1,(p+2)g_1]}}\right ).
\end{array}
$$
By Lemma~\ref{nuevo} and Remark~\ref{sumas} it is enough to verify that
the term:
\begin{equation}\label{aux}
\sum \limits_{t=0}^{1} c_t w_1 p^{3-t} v^{(t+p)g_1+2 } [1,(p+2-t) g_1 - 1]
+ c_{k-1} w_p p v^{(p+2)g_1 + 2 } [1, pg_1 -1],
\end{equation}
is an element of $\partial_3 (\low)$. The term
$\partial_3 (c_0 w_1 v^{pg_1 + 2} [1, (p+2)g_1 - 1] )$
takes care of the case $t=0$ in the expression (\ref{aux}).
When $t=1$ it is easily seen that the following equality holds:
$$
p^2 v^{g_2 + 2} [1,g_2 - 1] = \sum \limits_{i=0}^{pg_1} x_i p v^{g_2+g_1+2+i} [1,pg_1 - (i+1)].
$$
Therefore we have to prove that:
$$
p v^{g_2+g_1+2} [1,pg_1-1] =0 \;\;{\rm mod}\;\; \partial_3(\low),
$$
but an easy inductive argument in the second factor gives the proof.
For the term $u_1^{-1} p v^{pg_1} B^{[p,g_2+1]}$ we have that:
$$
u_1^{-1} \sum \limits_{j=1}^{pg_1} w_{j+g_1} p^2 v^{g_2 + j} [p, pg_1+1-j] =
u_1^{-1} \sum \limits_{t=0}^{g_{{}_1}^2 } x_t p v^{(p+2) g_1 + t}  [p, g^2_1 - t]
$$
and by Lemma~\ref{nuevo} this last expression is zero modulo $\partial_3(\low)$.
Therefore we have for $k=2$ that $\partial_{k+1} (p[p,(pk+1)g_1 +1] + \low)$ is equal to:
$$
y_1 p^3 v^{g_1} [1, (2p+1)g_1 + 1] +
u_1 y_1 p^2 v^{2 g_1} [1, 2pg_1 + 1] -
u_1^{-1} y_1 p^2 v^{g_2} [1, g_2 + 1]
$$
and this last expression is equal to:
$$
u_1^{-1} y_1 v^{2g_2} [1,1]
- y_1 p v^{g_1} (A+B)^{[1,(2p+1)g_1 + 1]}
- u_1^{-1} y_1 p v^{g_2} (A+B)^{[1,g_2 + 1]}.
$$
We verify that the sum
$$
\sum \limits_{i=1}^{p-2} \partial_{k+1} (y_1 v^{g_1+i} [1,(2p+1)g_1 + 1-i]),
$$
takes care of the term $-y_1 p v^{g_1} A^{[1,(2p+1)g_1 + 1]}$.

On the other hand, each summand of $-y_1 p v^{g_1} B^{[1,(2p+1)g_1 + 1]}$
is a Smith morphism image of (\ref{aux2}). The same argument applies to
each summand of $- u_1^{-1} y_1 v^{g_2} A^{[1,g_2 + 1]}$.

Finally, note that each term in the sum $- u_1^{-1} y_1 v^{g_2} B^{[1,g_2 + 1]}$
is a Smith morphism image of $p v^{g_2+g_1+1} [1,pg_1]$ and we have proved
that this term is zero modulo $\partial_3(\low)$. The theorem follows for
$k=2$. It is not difficult to verify that the theorem can be proved
for low values of $k$ using the same techniques.
Now suppose that $k \geq 6$; by Remark~\ref{bla} we have that $\partial_{k+1} (p [p, (pk+1)g_1 +1] + \low ) $ is
equal to:
$$
  y_1v^{g_1} \left (
{\cal S}_{k-1}^{[1,(pk+1)g_1+1]} \; - \; p^k {\cal S}_{k-1,k-6}^{[1,(pk+1)g_1 +1]}
\; + \;  p^k p_{k-5} (0) u_1^{-(2k-5)} v^{f_{k-5} g_1} [1, kg_1+ 2g_1^2]\right ) .
$$
It is easy to verify that modulo $\partial_{k+1}(\low)$ this term is equal to:
$$
\begin{array}{rl}
y_1 v^{g_1} & \left \{ (-1)^{k+1} (p_{k-5}^k(0) c_k(1) - p_{k-5}^{k-1} (1) )
u_1^{-k+2} p v^{(kp-g_1) g_1} [1, g_1^2 + g_1 + 1] \right . \\
\\
            & \left . (-1)^{k+1} ( p_{k-5}^k(0) d_k(1) - p_{k-5}^{k-2} (1) ) u_1^{-k+1} v^{(kp+1) g_1} [1,1] \right \}.
\end{array}
$$
Since the coefficient $a_{i,i+2} = 1$ for any $i$, the following holds:
$$
\begin{array}{lrl}
 & p_{k-5}^k(0) c_k(1) - p_{k-5}^{k-1} (1) & = \sum \limits_{i=1}^{k-3} a_{k-5,i}
 \left ( \sum \limits_{t=1}^{i+1}  d_{k-t} (0) \right ) \\
\\
 & & = \sum \limits_{i=1}^{k-4} (i+1) a_{k-5,i} + (k-2). \\
\\
 & p_{k-5}^k(0) d_k(1) - p_{k-5}^{k-2} (1) & = \sum \limits_{i=1}^{k-4} a_{k-5,i} (i+1) a_{k-5,i} + (k-3).
\end{array}
$$
Adding a suitable multiple of $\partial_{k+1} ([1, g_1^2 + (k+1)g_1 +1])$, we obtain:
$$
\partial_{k+1} (p [p,(pk+1)g_1 + 1] + \low) = (-1)^{k+1} y_1 u_1^{-k+1} v^{(kp+2)g_1} [1,1].
$$
\QED

\begin{Propo} \label{CP}
We have that $\partial_{k+1}([a,kg_1] + \low) = 0 $.
\end{Propo}

\section{Sharper relations in $F \otimes_{ku_*} ku_*(\mz_{p^2})$.}\label{s6}

In order to prove Proposition~\ref{CP} we need a sharper description of the elements of
$p^k$-torsion in $F \otimes_{ku_*} ku_*(\mz_{p^2})$.
This section is focused in the construction of certain polynomials that will be useful
in the combinatorics of the
non filtered differentials of the spectral sequence we are dealing with.
Some of this polynomials were used in \cite{kuosa}.
We begin this section establishing an easy result
about the behavior of ``{\it low }'' filtration terms with ``{\it high}'' $p$-divisibility.
We analyze the terms with $p^2$-divisibility in order to construct a stronger version of Theorem~\ref{Teoosa}.
For this purposes we will define a new family of polynomials. In this section
for $0 \leq i \leq p^2 - 1$ we will denote by $a_i = w_i p^{t_i} v^i$ the $i$-th coefficient of the
$[p^2]$-series. Here $w_i \in ku_*$ is a unit; as in the other sections $u_1$ denotes the unit $w_{p-1}$.

\bigskip

For $0 \leq k \leq g_2$ we define $q_k (x_1, \ldots, x_k) \in \mz_{(p)} [x_1, \ldots, x_k]$. We set $q_0 = -1$
and:
$$
q_k(x_1, \ldots, x_k) = \left \{ \begin{array}{lcl}
                                  - \sum \limits_{i=0}^{k-1} x_{k-i} q_i(x_1, \ldots, x_i) & {\rm for} & 1 \leq k \leq p-2 \\
                                  - \sum \limits_{i=0}^{p-2} x_{k-i} q_i(x_1, \ldots, x_i) & {\rm for} & p-1 \leq k \leq p^2-1
                                 \end{array}
\right .
$$

For $0 \leq i \leq p-2$ we will denote by:
$$
q_i^{[0]} = q_i (w_1, w_2, \ldots, w_i).
$$

Suppose $3 \leq n \leq p+1$. Let $0 \leq t \leq p-1$ and $0 \leq i \leq p-2$. We define:
$$
\widehat{q} \, {}^{\bf [0]}_{{\bf (t+1) g_1} + i} = \sum \limits_{k=0}^{p-2} q_k^{[0]}
w_{(t+1) g_1 + i - k} - q^{{}^{[0]}}_{{}_{\overline{i-t}}} w_{(t+1)p-1}.
$$
Here the polynomial
$ q^{{}^{[0]}}_{{}_{\overline{i-t}}} = q_k^{[0]}$ if $(i-t) \equiv k \,\, {\rm mod} \,\, p$
with $0 \leq k \leq p-2$ and it is zero if $(i-t) \equiv p-1 \,\, {\rm mod} \,\, p$. \\

For $0 \leq i \leq p-2$ and $0 \leq t \leq n-3$ we define:
$$
\widehat{q} \, {}^{\bf [1]}_{{\bf (t+1) g_1} + i} = \sum \limits_{k=0}^{i} q_{i-k}^{[0]} \,\,
\widehat{q} \, {}^{\bf [0]}_{{\bf (t+1) g_1} + k}.
$$
And for $0 \leq i \leq p-3$ we define:
$$
\widehat{q} \, {}^{\bf [1]}_{{\bf (t+1) g_1} + g_1+i} = \sum \limits_{k=i+1}^{p-2} q_{k}^{[0]} \,\,
\widehat{q} \, {}^{\bf [0]}_{{\bf (t+1) g_1} + (g_1+i-k)}.
$$

For $2 \leq s \leq n-2$ we denote by $H_{s,r}$ the set of $s$-partitions of $r$, with $s \leq r \leq n-2$.
For an element $ \iota \in H_{s,r}$ with $\iota = (k_1, \ldots, k_s)$ and for $0 \leq i \leq s(p-2)$
we define:
$$
\widehat{q} \, {}^{\bf [s]}_{\iota + i} = \sum \limits_{t=0}^{i'}
\widehat{q} \, {}^{\bf [0]}_{{\bf k_1 g_1} + t} \,\,
\widehat{q} \, {}^{\bf [s-1]}_{\iota'+i-t},
$$
where $i' = \min(i, p-2)$.

If we denote by $s' = s(p-2) + 1$ and take $0 \leq i \leq p-3$ we define the polynomial:
$$
\widehat{q} \, {}^{\bf [s]}_{\iota + s' + i} = \sum \limits_{k=i+1}^{p-2}
\widehat{q} \, {}^{\bf [0]}_{{\bf k_1 g_1} + k} \,\,
\widehat{q} \, {}^{\bf [s-1]}_{\iota'+ s' + i-k}.
$$
Here $\iota' \in H_{s-1,r-k_1}$ and is given by $\iota' = (k_2, \ldots, k_s)$. The polynomials
$\widehat{q} \, {}^{\bf [s-1]}_{\iota'+k}$ are defined by induction.

We define for $0 \leq i \leq p-3$ the polynomial:
$$
\widehat{q} \, {}^{\bf [0]}_{{\bf g_2} + i} = \sum \limits_{k=i+1}^{p-2} q_k^{[0]} w_{g_2 - k + i}.
$$
Recall that $g_2 = p^2 -1$.

The unique element of set $H_{n,n}$ will be denoted by $\iota_n = (1, \ldots, 1)$.

\bigskip

A direct calculation proves that for $n \geq 1$, the polynomial $q_0^{[n]} = (-1)^{n+1}$. We can
use an easy inductive argument to prove the following result.

\begin{Lema}
For $k \geq 3$ we have:
$$
p^k [a,kg_1] = \sum \limits_{i=0}^{p-2} u_1^{k-1} q_i^{[k-2]} p v^{(k-1) g_1 + i} [a, g_1 -i].
$$
in the tensor product $F \otimes_{ku_*} ku_* ({\mz}_{p^2})$.
\end{Lema}

\medskip

\begin{Teo}\label{part}
For $3 \leq n \leq p+1$ the term $p^2[a, ng_1] \in F \otimes_{ku_*} ku_*(\mz_{p^2})$ is equal to:
$$
\begin{array}{rl}
& \sum \limits_{k=1}^{n-1}  \sum \limits_{i=0}^{p-2} w_{kp - 1} q_i^{[0]} p v^{kp + i -1} [a, (n-k)g_1 - (i+k-1)] \\
+& \sum \limits_{j=1}^{n-2}  \sum \limits_{t=1}^{j} \sum H_{n-j-1,n-t-1} [a, (n-k)g_1 - \iota - (i+k-1)].
\end{array}
$$
Here the term $\sum H_{n-j-1,n-t-1}$ denotes
the sum:
$$
\sum \limits_{k=1}^{t}  \sum \limits_{i=0}^{n_j}  \sum \limits_{\iota }
w_{kp - 1} \,\,  \widehat{q} \, {}^{\bf [n-j-1]}_{\iota + i} p v^{\iota + kp +i -1},
$$
where $n_j = (n-j)(p-2)$ and $\iota \in H_{n-j-1,n-t-1}$.
\end{Teo}

\dem We proceed by induction on $n$.
It is not difficult to verify that by the relation imposed in the second factor
we have the following equality:
\begin{equation}\label{n3}
\begin{array}{rl}
p^2[a,3g_1] & = \sum \limits_{i=0}^{p-2} u_1 q_i^{[0]} p v^{g_1+i} [a, 2 g_1-i] \\
& + \sum \limits_{i=0}^{p-2} \widehat{q} \, {}^{\bf [0]}_{{\bf g_1} + i} p^2 v^{g_1+i} [a, 2g_1-i] \\
& + \sum \limits_{i=0}^{p-2} u_1 q_i^{[0]} w_{2p-1} p v^{2p-1+i} [a, g_1 - (i+1)]
\end{array}
\end{equation}
We know by Lemma 5.2 of \cite{kuosa} that:
$$
p^2 [a,2g_1] = \sum \limits_{i=0}^{p-2} u_1 q_i^{[0]} p v^{g_1+i} [a, g_1-i].
$$
Using this expression we obtain that:
$$
\widehat{q} \, {}^{\bf [0]}_{{\bf g_1} + i} p^2 v^{g_1+i} [a, 2g_1-i] =
 \sum \limits_{k=0}^{p-2} u_1 q_k^{[0]} \widehat{q} \, {}^{\bf [0]}_{{\bf g_1} + i} p v^{2g_1+k+i} [a, g_1-(i+k)].
$$
Therefore the second row of the expression (\ref{n3}) is equal to:
$$
\sum \limits_{i=0}^{p-2} \sum \limits_{k=0}^{i} u_1 q_{i-k}^{[0]}
\widehat{q} \, {}^{\bf [0]}_{{\bf g_1} + k} p v^{2g_1+i} [a, g_1-i].
$$
By the definition of the polynomials $\widehat{q} \, {}^{\bf [0]}_{{\bf g_1} + k}$ we have that:
$$
\begin{array}{rl}
p^2[a,3g_1] & =  \sum \limits_{k=0}^{2}
\sum \limits_{i=0}^{p-2} w_{kp-1} q_i^{[0]} p v^{kp+i-1} [a, (3-k)g_1-(i + k - 1)] \\
& + \sum \limits_{i=0}^{p-2} \widehat{q} \, {}^{\bf [1]}_{{\bf g_1} + i} p v^{2g_1+i} [a, g_1-i]
\end{array}
$$
Since $w_{p-1} = u_1$ and $H_{1,1} = \{ (1) \}$ the result follows for $n=3$. We suppose the result valid for
$3 \leq n < p+1$. We have that $p^2[a, (n+1)g_1] $ is equal to:
$$
\begin{array}{rl}
&\sum \limits_{k=0}^{n}
\sum \limits_{i=0}^{p-2} w_{kp-1} q_i^{[0]} p v^{kg_1+i+(k-1)} [a, (n-k+1)g_1-(i+k-1)] \\
&+\sum \limits_{i=0}^{p-2} \widehat{q} \, {}^{\bf [0]}_{{\bf g_1} + i} p^2 v^{g_1+i} [a, ng_1-i]
+ \sum \limits_{i=0}^{p-2} \widehat{q} \, {}^{\bf [0]}_{{\bf2 g_1} + i} p^2 v^{2g_1+i} [a, (n-1)g_1-i] \\
& + \ldots \\
&+ \sum \limits_{i=0}^{p-2}  \widehat{q} \, {}^{\bf [0]}_{{\bf (n-2) g_1} + i} p^2 v^{(n-2)g_1+i} [a, 3g_1-i]
+ \sum \limits_{i=0}^{p-2}  \widehat{q} \, {}^{\bf [0]}_{{\bf (n-1) g_1} + i} p^2 v^{(n-1)g_1+i} [a, 2g_1-i]
\end{array}
$$
By Lemma 5.2 of~\cite{kuosa} we have that $p^2[a,g_1] =0$. Also note that depending on the value of $p$ some terms can be zero.
Therefore we have for $1 \leq k \leq n-2$:
$$
\sum \limits_{i_0=0}^{p-2}  \widehat{q} \, {}^{\bf [0]}_{{\bf k g_1} + i_0} p^2 v^{kg_1+i_0} [a, (n-k+1)g_1-i_0]
$$
is equal to:
$$
\begin{array}{rl}
&\sum \limits_{i_0=0}^{p-2}  \widehat{q} \, {}^{\bf [0]}_{{\bf k g_1} + i_0} v^{kg_1+i_0}
\sum \limits_{t=1}^{n-k}  \sum \limits_{i_1=0}^{p-2} w_{tp - 1} q_{i_1}^{[0]}
p v^{tp + i_1 -1} [a, (n-k-t+1)g_1 - (i_0+i_1+t-1)] \\
+&  \sum \limits_{i_0=0}^{p-2}  \widehat{q} \, {}^{\bf [0]}_{{\bf k g_1} + i_0} v^{kg_1+i_0}
\sum \limits_{j=1}^{n-k-1}  \sum \limits_{t=1}^{j} \sum H_{n-k-j,n-k-t}
[a, (n-k_0)g_1 - \iota - (i_0+i+k_0-1)].
\end{array}
$$
Here $1 \leq k_0 \leq t$ is the parameter that arise from the sum $H_{n-k-j,n-k-t}$.

On the other hand we have the following:
$$
\sum \limits_{i_0=0}^{p-2} \sum \limits_{i_1=0}^{p-2}
\widehat{q} \, {}^{\bf [0]}_{{\bf k g_1} + i_0}  q_{i_1}^{[0]}
= \sum \limits_{i=0}^{2(p-2)} \sum \limits_{i_0+i_1=i}
\widehat{q} \, {}^{\bf [0]}_{{\bf k g_1} + i_0} q_{i_1}^{[0]}
=\sum \limits_{i=0}^{2(p-2)} \widehat{q} \, {}^{\bf [1]}_{{\bf k g_1} + i}.
$$

For fixed natural numbers $1 \leq j \leq n-k-1$ and $1 \leq t \leq j$ we have:
$$
\sum \limits_{i_0=0}^{p-2} \widehat{q} \, {}^{\bf [0]}_{{\bf k g_1} + i_0}
\sum H_{n-k-j,n-k-t} =
\sum \limits_{i_0=0}^{p-2}
\sum \limits_{k'=1}^{t}  \sum \limits_{i_1=0}^{(n+1)_{{}_{k+j}}}
\sum \limits_{\iota } w_{k'p - 1} \,\,
\widehat{q} \, {}^{\bf [0]}_{{\bf k g_1} + i_0}
\widehat{q} \, {}^{\bf [n-k-j]}_{\iota + i_1} p v^{\iota + k'p +i_1 -1}.
$$
Here $\iota \in H_{n-k-j,n-k-t}$. Therefore, for $1 \leq k' \leq t$ we have:
$$
\sum \limits_{i_0=0}^{p-2} \sum \limits_{i_1=0}^{(n+1)_{{}_{k+j}}}\sum \limits_{\iota }
\widehat{q} \, {}^{\bf [0]}_{{\bf k g_1} + i_0} \,\,
\widehat{q} \, {}^{\bf [n-k-j]}_{\iota + i_1}
= \sum \limits_{\iota} \sum \limits_{i=0}^{(n+2)_{{}_{k+j}}}
\sum \limits_{i_0+i_1=i} \widehat{q} \, {}^{\bf [0]}_{{\bf k g_1} + i_0} \,\,
\widehat{q} \, {}^{\bf [n-k-j]}_{\iota + i_1}.
$$
Since $0 \leq i_0 \leq p-2$ and $0 \leq i_1 \leq (n+1)_{k+j}$,
for a fixed $\iota \in H_{n-k-j,n-k-t}$:
$$
\sum \limits_{i=0}^{(n+2)_{{}_{k+j}}} \sum \limits_{i_0+i_1=i}
\widehat{q} \, {}^{\bf [0]}_{{\bf k g_1} + i_0} \,\,
\widehat{q} \, {}^{\bf [n-k-j]}_{\iota + i_1} =
\widehat{q} \, {}^{\bf [n-k-j+1]}_{\iota_1 + i}.
$$
Here $\iota_1 = (k,\iota) \in H_{n-k-j+1,n-t}$ and $0 \leq i \leq (n+2)_{k+j}$.
\smallskip

Note that for a given element $\kappa \in H_{n-j,n-t}$
with $1 \leq j \leq n-2$ and $1 \leq t \leq j$, if $\kappa = (a_1, a_2, \ldots, a_{n-j})$, then
$(a_2, \ldots, a_{n-j}) \in  H_{n-j-1,n-t-a_1}$, and since
$$
1 \leq a_1 \leq n-t -(n-j-1) = j-t+1
$$
there exist $1 \leq h \leq n-2$ such that:
$$
\kappa = (k, \kappa') \,\,\,\, {\rm with} \,\,\,\, \kappa' \in H_{n-j-1,n-t-a_1}.
$$
Finally note that $
\sum \limits_{i=0}^{p-2}
\widehat{q} \, {}^{\bf [0]}_{{\bf (n-1) g_1} + i} \,\, p^2 v^{(n-1)g_1 +i} [a,2g_1 - i]$ is equal to:
$$
\sum \limits_{i_0=0}^{p-2} \sum \limits_{i_1=0}^{p-2} u_1
 \widehat{q} \, {}^{\bf [0]}_{{\bf (n-1) g_1} + i_0} \, q_{i_1}^{[0]} p v^{ng_1 +i_0 + i_1} [a, g_1 -(i_0+i_1)] +
\sum \limits_{i=0}^{p-2} u_1 \widehat{q} \, {}^{\bf [1]}_{{\bf (n-1) g_1} + i}
 p v^{ng_1 + i} [a,g_1-i].
$$
This completes the proof of the theorem, since the set $H_{n-1,n-1}$ has one element.
\QED

Now we analyze the $p^2$-torsion elements of arbitrary filtration.
For natural numbers $s$ and $r \geq p+2$ we denote by $H'_{s,r}$ the set of $s$-partitions of $r$, whose
entries are less or equal to $p+1$.
For an element $\iota \in H'_{s,r}$ and $0 \leq i \leq (s+1)(p-2)$ we can define the polynomial $\widehat{q} \, {}^{\bf [s]}_{\iota + i}$,
using the polynomials constructed in the beginning of the section.

An easy consequence of Theorem~\ref{part} is the following result.

\begin{Coro}\label{relCP}
For $n \geq p+2$ the element $p^2[a,ng_1]$ is equal to:
$$
\begin{array}{rl}
& \sum \limits_{k=1}^{p-1}  \sum \limits_{i=0}^{p-2} w_{kp - 1} q_i^{[0]} p v^{kp + i -1} [a, (n-k)g_1 - (i+k-1)] \\
+& \sum \limits_{j=1}^{n-2}  \sum \limits_{t=1}^{j} \sum H'_{n-j-1,n-t-1} [a, (n-k)g_1 - \iota - (i+k-1)]\\
+& \sum \limits_{i=0}^{p-2} u_2 q_i^{[0]} v^{p^2 + i -1} [a, ng_1 - g_2 -i] \\
+& \sum \limits_{j=1}^{n-p-2}  \sum \limits_{t=1}^{j} \sum
H^{''}_{n-p-j-1,n-p-t-1} [a, ng_1 - g_2 - \iota - i]
\end{array}
$$
Here the term $\sum H'_{n-j-1,n-t-1}$ denotes the sum:
$$
\sum \limits_{k=1}^{t}  \sum \limits_{i=0}^{n_j}  \sum \limits_{\iota }
w_{kp - 1} \,\,  \widehat{q} \, {}^{\bf [n-j-1]}_{\iota + i} p v^{\iota + kp +i -1},
$$
where $n_j = (n-j)(p-2)$ and $\iota \in H'_{n-j-1,n-t-1}$.
The term $\sum H^{''}_{n-p-j-1,n-p-t-1}$ is equal to:
$$
\sum \limits_{i=0}^{n_j}  \sum \limits_{\iota } u_2 \,\,  \widehat{q} \, {}^{\bf [n-j-1]}_{\iota + i}  v^{\iota + p^2 +i -1},
$$
where $n_j = (n-j)(p-2)$ and $\iota \in H'_{n-p-j-1,n-p-t-1}$.
\end{Coro}

\dem (Proof of Proposition~\ref{CP}) Note that we can take $a = n g_1$ for some natural number $n$. Let us consider the following sum:
$$
\begin{array}{l}
\,\,\,\, \sum \limits_{i=0}^{n-1} \partial_{k+1} \left (  u_1^{-i} y_1^{i} [(n-i)g_1 , (k+i) g_1 ] \right ) \\
+ \sum \limits_{i=0}^{n-1}
\partial_{k+1} \left ( u_1^{-(i+2)} y_1^{i}  y_2  v^{g_2-2g_1} [(n-i)g_1 - g_2, (k+i+2) g_1 ] \right ) \\
\vdots \\
+ \sum \limits_{i=0}^{n-1}
\partial_{k+1} \left ( u_1^{-(i+k)} y_1^{i} y_k v^{g_k - kg_1} [(n-i)g_1-g_k , (2k+i) g_1 ] \right ) + \partial_{k+1} (\bf{SM}).
\end{array}
$$
Here ${\bf SM}$ denotes elements that are Smith morphism images of the first rows. Using Corollary~\ref{relCP} it
is not difficult to verify that this sum is equal to:
$$
\begin{array}{rl}
& u_1^{-1} u_2 y_k \sum \limits_{i=0}^{p-2} q_i^{[0]} v^{g_k +g_2 - kg_1 +i} [ng_1 - g_k , (k+1) g_1 - g_2 -i] \\
+ &u_1^{-2} u_2 y_1 y_k \sum \limits_{i=0}^{p-2} q_i^{[0]} v^{g_k +g_2 - kg_1 +i} [(n-1)g_1 - g_k , (k+2) g_1 - g_2 -i] +
\cdots + \\
+ & u_1^{-1} u_2 y_k \sum \limits_{i=0}^{p-2} q_i^{[0]} v^{g_k +g_2 - kg_1 +i} [4g_1 - g_k , (k+1) g_1 - g_2 -i] \\
+ & y_{k+1} v^{g_{k+1}} [n g_1 - g_{k+1} , k g_1] + u_1^{-1} y_1 y_{k+1} v^{g_{k+1}} [(n-1) g_1 - g_{k+1} , (k+1) g_1] + \cdots + {\bf SM}.
\end{array}
$$
Now we will prove that each of these terms are image under $\partial_{k+1}$ of ``$\low$'' filtration terms.

\sc

\bigskip\medskip
Leticia Z\'arate

CEFyMAP - Universidad Aut\'onoma de Chiapas.

4a. Oriente Norte No. 1428.  Entre 13a. y 14a. Norte.

Col. Barrio La Pimienta

Tuxtla Guti\'errez, Chiapas. C.P. 29000

M\'exico

e-mail: {\tt leticia@math.cinvestav.mx}


\begin{thebibliography}{rrr}

\bibitem{Nakos}{G.~Nakos: On the Brown-Peterson
homology of certain classifying spaces,
Ph.D. Thesis, The Johns Hopkins University, 1985.}

\bibitem{tesis}{L.~Z\'arate: On the $BP \langle n \rangle_*$-homology
of $\mz_{2^e} \times \mz_{2^e}$,
Ph.D. Thesis, CINVESTAV - IPN, 2007.}

\bibitem{BP-osa}{L.~Z\'arate: {\it On the $BP$-homology
of $\mz_{2^e} \times \mz_{2^e}$},
Journal of K-theory, (To appear).}

\bibitem{kuosa}{L.~Z\'arate: {\it On the $ku$-homology
of certain classifying spaces}, (Submitted).}

\end{thebibliography}
\end{document}